\documentclass[12pt]{article}
\usepackage{amssymb,oldlfont,upref}

\title{Regularized traces and Taylor expansions for the heat semigroup}
\author{Michael Hitrik\\Department of Mathematics\\University of California\\Berkeley, 
CA 94720\\hitrik@math.berkeley.edu\and 
Iosif Polterovich\\MSRI\\1000 Centennial Dr.\\Berkeley, CA 94720\\iossif@msri.org}
\def\wrtext#1{\relax\ifmmode{\leavevmode\hbox{#1}}\else{#1}\fi}
\def\abs#1{\left|#1\right|}
\def\begeq{\begin{equation}}
\def\endeq{\end{equation}}
\def\Remark{\vskip 2mm \noindent {\em Remark}}

\newcommand{\eps}{\varepsilon}
\def\part#1{\frac{\partial}{\partial #1}}
\def\half{\frac{1}{2}}

\def\norm#1{||\,#1\,||}

\newcommand{\real}{\mbox{\bf R}}

\newcommand{\nat}{\mbox{\bf N}}

\newcommand{\Spec}{\mbox{\rm Spec\,}}
\newcommand{\ad}{\mbox{\rm ad\,}}

\newcommand{\Tr}{\mbox{\rm Tr\,}}
\renewcommand{\exp}{\mbox{\rm exp\,}}

\newtheorem{dref}{Definition}[section]
\newtheorem{lemma}[dref]{Lemma}
\newtheorem{theo}[dref]{Theorem}
\newtheorem{prop}[dref]{Proposition}

\newenvironment{proof}{\vspace{.3cm}\noindent{{\em Proof:}}}{\hfill$\Box$}
\begin{document}
\maketitle
\begin{abstract}
We compute the coefficients in asymptotics of regularized traces and associated trace (spectral) distributions for 
Schr\"odinger operators, with short and long range potentials. A kernel expansion for the Schr\"odinger semigroup is derived, and 
a connection with non-commutative Taylor formulas is established. 
\end{abstract}
\section{Introduction and main results}
\setcounter{equation}{0}
\setcounter{dref}{0}
The purpose of this paper is to describe asymptotics of regularized traces and trace distributions associated to the 
Schr\"odinger operator $H=H_V=H_0+V$ in $\real^n$. Here $H_0=-\Delta$ and $V$ is a real-valued, smooth, and bounded 
potential satisfying some suitable decay conditions at infinity. In particular, in the short range case, we present closed 
formulas for the heat invariants of 
$H$, which then leads to formulas for the coefficients in the asymptotic expansion of the derivative of the scattering phase of 
$H$. 

Most of the paper is devoted to the long range case. Here we use the results of 
Melin~\cite{Melin} as our starting point. In~\cite{Melin} 
regularized 
traces of functions of $H$ were studied. We sharpen the main results 
of~\cite{Melin}, and then we calculate the coefficients arising in 
the asymptotic expansion of these traces and in the expansions of 
associated trace distributions. This is done by
applying the techniques for computing the local heat invariants of the Laplace-Beltrami operator, which were 
developed in~\cite{Pol1},~\cite{Pol3}. The latter works were
based on the resolvent kernel 
expansion of Agmon and Kannai~\cite{AK} derived by means of  some commutator 
formulas. A similar commutator technique for $H$ was then used in~\cite{Melin}. Here we eliminate multiple 
commutators of $H_0$ and $V$ from the expansions 
of~\cite{Melin}, just as in~\cite{Pol1} similar commutators were eliminated from the expansions of~\cite{AK}. This makes it possible 
to compute the coefficients  in the asymptotics of trace distributions. 
It also
provides an analogy between the expansions considered in~\cite{AK},
\cite{Melin} and Taylor series for functions of (Schr\"odinger) operators.

Our main results are as follows. First, we derive an asymptotic
expansion for the heat kernel on the diagonal, $e^{-tH}(x,x)$. It may be viewed as a heat version of the Agmon-Kannai expansion. We remark 
that a relation to Taylor expansions can be already seen in this result.  
\begin{theo}
\label{hh}
For every positive integer $N$, the following asymptotic representation of the heat kernel is true as $t\to 0^+$, 
locally uniformly in $x\in \real^n$, 
\begin{equation}
\label{great}
e^{-tH}(x,x)=\sum_{m=0}^{N} \frac{t^m}{m!} \left(X_m e^{-tH_0}\right)(x,x)+{\cal O}(t^{(N+2)/2-n/2}).
\end{equation}
Here the differential operators $X_m$ are defined by 
\begeq
\label{xj}
X_m=\sum_{k=0}^m (-1)^k {m \choose k} H^k H_0^{m-k}.
\endeq
When $N$ is odd, the remainder term in (\ref{great}) is ${\cal O}(t^{(N+3)/2-n/2})$. 
\end{theo}
We use Theorem \ref{hh} for studying asymptotics of regularized traces
for Schr\"odinger operators with long-range potentials.
Let $V\in S^{-\eps}(\real^n)$, for some $\varepsilon\in(0,1]$, so that  
\begeq
\label{VV}
\abs{\partial^{\nu} V(x)}\leq C_{\nu} (1+\abs{x})^{-\eps-\abs{\nu}},\quad x\in \real^n. 
\endeq
We also 
set $N:=[\frac{n}{\varepsilon}]$, the integer part of $n/{\varepsilon}$.
\begin{theo}
\label{regt}
The operator 
$$
e^{-tH}-\sum_{m=0}^{N}\frac{t^m}{m!}X_m e^{-tH_0}, \quad t>0, 
$$
is of trace class. As $t\to 0^+$, one has the asymptotic expansion 
\begeq
\label{reg}
\Tr\left(e^{-tH}-\sum_{m=0}^{N}\frac{t^m}{m!}X_m e^{-tH_0}\right)\sim (4\pi t)^{-n/2} \sum_{j=1}^{\infty} 
\alpha_j t^{j},
\endeq
where the coefficients $\alpha_j=\int \alpha_j(x)\,dx$, and the smooth functions $\alpha_j(x)\in L^1(\real^n)$ are defined as follows: 
\begeq
\label{alpha1}
\alpha_j(x)=0,\quad j<(N+2)/2,
\endeq
\begeq
\label{alpha2}
\alpha_j(x)=a_j(x):=
(-1)^j\sum_{k=0}^{j-1}
{j-1+\frac{n}{2}\choose k+\frac{n}{2}} 
\frac{H_{y}^{k+j}\left(d(x,y)^{2k}\right)|_{y=x}}{4^k k! (k+j)!},\quad
j\ge N+1,
\endeq
where $d(x,y)$ is the Euclidean distance between $x$, $y\in \real^n$. When $(N+2)/2\le j \leq N$ we have
\begeq
\label{alpha3}
\alpha_j(x)=a_j(x)-
(-1)^j\sum_{k=0}^{N-j}
{N-j+\frac{n}{2}\choose k+\frac{n}{2}} 
\frac{H_{y}^{k+j}\left(d(x,y)^{2k}\right)|_{y=x}}{4^k k! (k+j)!}
\endeq
\end{theo}
{\em Remark.} The coefficients $a_j(x)$ defined in (\ref{alpha2}) are the local heat invariants of $H$ --- see Theorem 3.1. 

The plan of the paper is as follows. In Section 2 we prove Theorem \ref{hh}.
The formulas for the scattering phase coefficients in the short range case 
are presented in Section 3, and Section 4 is then devoted to the case of long range potentials. Here we first improve some of 
the results of~\cite{Melin}. These improvements  are then used to 
prove Theorem \ref{regt} and to describe asymptotics of 
the associated trace distributions. In Section 5 we relate the heat kernel expansions of Theorem 1.1 to the Taylor 
expansion for the heat semigroup and to a non-commutative Taylor formula of~\cite{Kantor1}. 
Finally, Section 6 is devoted to the proof of an auxiliary 
technical proposition which is instrumental in establishing the results of Section 4.

\noindent{\it Acknowledgements.} We are very grateful to Yakar Kannai and Anders Melin for their interest in  this work. 
We would also like to thank John Lott, Kate Okikiolu, Alexander Pushnitski, 
and Barry Simon for useful discussions. 
The second author is indebted to Victor and 
Leonid Polterovich for helpful advice. 
The first author is supported by the Swedish Foundation for International
Cooperation in Research and Higher Education (STINT). The second author is supported by the CRM/ISM (Montr\'eal) and MSRI (Berkeley) 
postdoctoral fellowships. 

\section{Heat kernel expansions}
\setcounter{equation}{0}
\setcounter{dref}{0}

In this section we are going to consider the Schr\"odinger operator $H=H_0+V$ in $\real^n$, where $H_0=-\Delta$ 
and $n\geq 1$ is arbitrary. Here the real-valued potential $V$ will be $C^{\infty}$-smooth, and for simplicity we shall assume that it is a 
bounded function on $\real^n$. Then $H$ is a self-adjoint operator on $L^2(\real^n)$ with the domain $D(H)=D(H_0)=H^2(\real^n)$, the 
standard Sobolev space on $\real^n$. 

We shall consider the heat semigroup $e^{-tH}$, $t\geq 0$, acting on $L^2(\real^n)$. It is well-known that under the assumptions above, the 
distribution kernel of $e^{-tH}$, which will be denoted $e^{-tH}(x,y)$, is a smooth function of $(x,y)$ which depends continuously on $t>0$. 
Moreover, for every integer $N\geq 1$, we have, as $t\rightarrow 0+$, 
\begeq
\label{eq1}
e^{-tH}(x,x)=(4 \pi t)^{-n/2}\left(\sum_{j=0}^{N} a_j(x) t^j+{\cal O}_x(t^{N+1})\right), 
\endeq
locally uniformly in $x\in \real^n$. The smooth functions $a_j(x)$ will be called the local heat invariants of $H$. 
We notice that the existence of the expansion (\ref{eq1}) can be inferred from the Feynman-Kac representation of 
$e^{-tH}(x,x)$, see~\cite{Simon}, 
\begeq
\label{eq2}
e^{-tH}(x,x)=(4 \pi t)^{-n/2} E\left(\exp\left(-\int_0^t V(x+\kappa(s))\,ds\right)\right),
\endeq
where $\kappa(s)$, $0\leq s\leq t$, is the Brownian bridge with $\kappa(0)=0$, which returns to $0$ at time $t$. We recall that it is a 
Gaussian process of mean zero, and covariance $E(\kappa_j(s)\kappa_k(u))=2s(1-u/t)\delta_{jk}$, $0\leq s\leq u\leq t$. Taking a Taylor 
expansion of the exponential in the right-hand side of (\ref{eq2}), and then expanding $V$ in a Taylor series around $x$ leads to (\ref{eq1}). 
(See ~\cite{SimonandCo} for further applications of this idea to some one-dimensional Schr\"odinger operators, and also, Section 6.) We get 
the following expressions for the first coefficients in (\ref{eq1}), 
$$
a_0(x)=1,\quad a_1(x)=-V(x), \quad a_2(x)=\frac{1}{6}\left(3 V^2(x)-\Delta V(x)\right). 
$$
Theorem \ref{hh} provides an asymptotic expansion for $e^{-tH}(x,x)$ well suitable for computation of all the coefficients $a_j(x)$. See also 
Section 3.

\vskip 2mm
\noindent {\em Proof of Theorem 1.1:} 
It follows easily from (\ref{xj}) that the operators $X_m$ satisfy
$$
X_m=(H_0-H)X_{m-1}+[X_{m-1},H_0],\quad m\geq 1.
$$
Since the order of the operator $H_0-H=-V$ is zero, a simple proof by induction shows that the order of the differential operator 
$X_m$ is $\leq m-1$, $m\geq 1$. We have that the sum in the right-hand side of (\ref{great}) is $t^{-n/2}$ times a polynomial of degree $N$ in 
$t$, and we shall now compute the coefficients there. Since  
$$
e^{-tH_0}(x,y)=\frac{1}{(4 \pi t)^{n/2}} e^{-\abs{x-y}^2/4t}, 
$$
it follows from Taylor's formula that when $\mu=(\mu_1,\ldots \mu_n)$ is an arbitrary multi-index, it is true that 
$$
\frac{\partial^{2\mu}}{\partial x^{2\mu}}e^{-tH_0}(x,x)=(4 \pi t)^{-n/2} \frac{(-1)^{\abs{\mu}} (2\mu)!}
{4^{\abs{\mu}} t^{\abs{\mu}} \mu!}.
$$
It follows that if we write 
$$
X_m=\sum_{\abs{\nu}\leq m-1} p_{\nu}(x) \partial^{\nu}_x,
$$
then 
$$
\left(X_m e^{-tH_0}\right)(x,x)= (4 \pi t)^{-n/2} \sum_{2\abs{\mu}\leq m-1}
p_{2\mu}(x) \frac{(-1)^{\abs{\mu}}(2\mu)!}{4^{\abs{\mu}} 
t^{\abs{\mu}} \mu!},
$$
since to get a non-zero contribution, all indices in $\nu=(\nu_1,\ldots \nu_n)$ must be even. Thus the total contribution of the term 
$\frac{t^m}{m!} X_m e^{-tH_0}(x,x)$ to the coefficient in front of $t^j$, $1\leq j\leq N$, is equal to 
$$
(4\pi t)^{-n/2}\frac{1}{m!} \sum_{\abs{\mu}=m-j} p_{2\mu}(x) \frac{(-1)^{\abs{\mu}}(2\mu)!}{4^{\abs{\mu}} \mu!}.
$$
In particular, $m\geq j$. On the other hand, since $2\abs{\mu}\leq m-1$ and $m\leq N$ it follows that $m\leq 2j-1$, provided that 
$j<(N+2)/2$. (If $N$ is odd this means that $j\leq (N+1)/2$.) We conclude that for these values of $j$, the coefficient in front of 
$t^j$ is given by $(4 \pi t)^{-n/2}$ times the sum 
\begeq
\label{eq3}
\sum_{m=j}^{2j-1} \frac{1}{m!} \sum_{\abs{\mu}=m-j} X_{m,y}\left(\frac{(y-x)^{2\mu}}{(2\mu)!}\right)\bigg |_{y=x}
\frac{(-1)^{\abs{\mu}}(2\mu)!}{4^{\abs{\mu}} \mu!}.
\endeq
Here we have also used that $p_{2\mu}(x)=X_{m,y}\left(\frac{(y-x)^{2\mu}}{(2\mu)!}\right)\bigg |_{y=x}$. The expression 
(\ref{eq3}) is precisely the formula for the $j$-th local heat invariant $a_j(x)$ of $H$, according to 
Theorem 3.1.1 of~\cite{Pol3}. We remark that while in~\cite{Pol3} such a formula was obtained for the heat invariants of the 
Laplace-Beltrami operator on a Riemannian manifold, it follows from the proof in~\cite{Pol3}, which is based on~\cite{AK}, that the 
$a_j(x)$ in our case are given by (\ref{eq3}). In particular, the upper summation index in the outer sum in (\ref{eq3}) is $2j-1$ instead of 
$4j$, as in~\cite{Pol3}, since in the case at hand, the order of the operator $X_m$ is $\leq m-1$. 
An application of (\ref{eq1}) completes the proof.  
{\hfill$\Box$}


\noindent
\Remark. The operators $X_m$ were introduced in~\cite{Pol1} in order to simplify expansions of the resolvent kernels of 
elliptic self-adjoint differential operators, due to Agmon and Kannai~\cite{AK}. In particular, as a consequence of the results of~\cite{AK}, 
it was proved in~\cite{Pol1} that when $l\in \nat$ is such that $l\geq n/2$, the following asymptotic expansion is true for the 
kernel of the $l$-th derivative of the resolvent $R(\lambda)=(H-\lambda I)^{-1}$, as $\lambda \rightarrow \infty$ outside a conic 
neighbourhood of the set $[E_0,\infty)$, where $E_0=\inf \Spec(H)$, 
\begeq
\label{resexp}
\left(\frac{d}{d\lambda}\right)^l R(\lambda)(x,x)\sim \sum_{m=0}^{\infty} \frac{(m+l)!}{m!} 
\left(X_m F(\lambda)^{m+l+1}\right)(x,x).
\endeq
Here $F(\lambda)=(H_0-\lambda I)^{-1}$. 
Using (\ref{resexp}), Theorem \ref{hh} can also be proved by means of the Cauchy's integral formula, expressing the heat semigroup 
in terms of $R(\lambda)$. Since in~\cite{Pol3} the local heat invariants were computed using 
(\ref{resexp}), the expansion (\ref{resexp}) is essentially 
equivalent to that of Theorem \ref{hh}. In our applications, we have found 
it more convenient to work with (\ref{great}).

\section{Asymptotics for short range potentials}
\setcounter{equation}{0}
\setcounter{dref}{0}

As a preparation for further considerations, in this section we shall work with compactly supported potentials. More precisely, we shall 
assume that the potential $V\in C_0^{\infty}(\real^n)$. It then follows easily from Duhamel's formula 
\begeq
\label{duhamel}
e^{-tH}-e^{-tH_0}=\int_0^t e^{-sH}V e^{-(t-s)H_0}\,ds,
\endeq
that the operator $e^{-tH}-e^{-tH_0}$ is of trace class, $t>0$. Moreover, it is well-known that the scattering matrix 
$S(\lambda): L^2(S^{n-1})\rightarrow L^2(S^{n-1})$, $\lambda>0$, of $H$ is a trace class perturbation of the identity, so that 
$\det S(\lambda)$ is well-defined. We may therefore introduce the scattering phase 
$$
s(\lambda)=\frac{1}{2 \pi i} \log \det S(\lambda), \quad \lambda>0,
$$
which is well-defined modulo the integers. The relation between the trace of (\ref{duhamel}) and $s(\lambda)$ 
is given by the Birman-Krein formula, see~\cite{Gp} and~\cite{Muller}, 
\begeq
\label{BKheat}
\Tr(e^{-tH}-e^{-tH_0})=\sum_{j=1}^l e^{-t\lambda_j}+\eps_0+\int_0^{\infty}e^{-t\lambda} s'(\lambda)\,d\lambda,\quad t>0. 
\endeq
Here $\lambda_j\leq 0$, $j=1,\ldots l$ are the eigenvalues of $H$. The constant $\eps_0$ comes from the 
contribution of the resonance state at the origin. 

An essentially well-known argument based on the Feynman-Kac formula (see~\cite{CdV}), together with (\ref{eq1}) shows that the trace 
$$
\Tr (e^{-tH}-e^{-tH_0})=\int \left(e^{-tH}(x,x)-e^{-tH_0}(x,x)\right)\,dx
$$
has an asymptotic expansion as $t\rightarrow 0^+$,
\begeq
\label{h1}
\Tr(e^{-tH}-e^{-tH_0})\sim (4\pi t)^{-n/2} \sum_{j=1}^{\infty} a_j t^j,\;\;t\rightarrow 0^+,
\quad a_j=\int a_j(x)\,dx. 
\endeq  
The coefficients $a_j$ will be called the heat invariants of $H$. We shall state here the main result of~\cite{Pol3}, adapted to the case of 
the Schr\"odinger operator $H$. (See also~\cite{Pol2} for computation of the heat coefficients in one dimension, where they are precisely the KdV 
invariants.)

\begin{theo}
\label{scat}
The heat invariants of $H$ are given by 
$$
a_j=\int_{{\small \bf R}^n} a_j(x)\,dx, \quad j\geq 1, 
$$
where
\begeq
\label{local}
a_j(x)=(-1)^j \sum_{k=0}^{j-1} {j-1+\frac{n}{2}\choose k+\frac{n}{2}} 
\frac{H_y^{k+j}\left(d(x,y)^{2k}\right)|_{y=x}}{4^k k! (k+j)!}. 
\endeq
Here $d(x,y)=|x-y|$ is the Euclidean distance between the points $x,y\in \real^n$, and the binomial coefficients for $n$ odd 
are given by 
$$
{j-1+\frac{n}{2} \choose k+\frac{n}{2}}=\frac{\Gamma(j+\frac{n}{2})}{(j-1-k)!\, \Gamma(k+\frac{n}{2}+1)}.
$$
\end{theo}

\Remark. We refer to~\cite{BsB} for a different set of expressions for the heat invariants of Schr\"odinger operators.

The asymptotic expansion (\ref{h1}) corresponds to an asymptotic expansion of the derivative of the scattering phase. We have that, 
depending on the parity of the dimension, 
\begeq
\label{odd}
s'(\lambda) \sim \sum_{j=1}^{\infty} b_j \lambda^{n/2-j-1},
\quad \lambda\rightarrow \infty,\;\;n\;\;\wrtext{odd}, 
\endeq
and
\begeq
\label{even}
s'(\lambda)=\sum_{j=1}^{n/2-1} b_j \lambda^{n/2-j-1}+{\cal O}(\lambda^{-\infty}), \quad \lambda\rightarrow 
\infty,\;\;n\;\;\wrtext{even}. 
\endeq
This result has been established in~\cite{BuslaevFaddeev},~\cite{Buslaev} in the one- and three-dimensional case, respectively, 
and then in more general situations in~\cite{CdV},~\cite{Gp},~\cite{Popov},~\cite{Robert}. The latter paper considers potentials that 
are standard symbols of order $-\rho$, $\rho>n$. 

An application of Lemma 5.2 of~\cite{CdV} to (\ref{BKheat}) shows that there is the following relation between 
the coefficients $a_j$ and $b_j$, 
\begeq
\label{relationab}
b_j=(4\pi) ^{-n/2} \frac{a_j}{\Gamma(n/2-j)}, \quad j=1,2,\ldots,\;\;n\;\wrtext{is odd},\;\;j=1,2,\ldots \frac{n}{2}-1,
\;\;n\;\wrtext{is even},
\endeq
and the $a_j$ are given in Theorem 3.1. 

\section{Trace distributions for long range potentials}
\setcounter{equation}{0}
\setcounter{dref}{0}
In this section we shall study asymptotics of trace distributions associated to the Schr\"odinger operator with 
long range potentials, introduced in~\cite{Melin}. The trace distributions generalize the scattering phase in the short range case. 
Following~\cite{Melin}, we shall assume that the potential $V(x)$ satisfies 
(\ref{VV}), i.e.
$$
\abs{\partial^{\nu} V(x)}\leq C_{\nu} (1+\abs{x})^{-\eps-\abs{\nu}},\quad x\in \real^n 
$$ 
for some $\varepsilon\in(0,1]$.
For this class of potentials (and, in fact, for a much larger class of pseudo-differential perturbations), expansions for functions 
of $H$ were obtained in~\cite{Melin}. In order to describe these expansions, we have to recall the definition of 
the operators $V_j$ from~\cite{Melin}. 

Let $J=(j_1,\ldots ,j_k)$ be a finite vector with nonnegative integer components. We set $\abs{J}=j_1+\ldots +j_k$ and 
$\norm{J}=\abs{J}+k$. An operator $V_J$ is inductively defined for all such vectors $J$ in the following way:
If $J$ is the empty sequence then $V_J=I$, the identity 
operator on $L^2$. Assuming that $J=(j_1,\ldots ,j_k)$ is non-empty and that 
$V_J$ has already been defined for all shorter sequences, we set
$$
V_J=\left(\ad H_0\right)^{j_k}\left(V_{(j_1,\ldots ,j_{k-1})}V\right),
$$
where $V_{(j_1,\ldots ,j_{k-1})}=I$ if $k=1$, and $(\ad H_0)A=[H_0,A]=H_0A-AH_0$. The operators $V_j$ are then defined by the formula 
\begeq
\label{vj}
V_j=\sum_{\norm{J}=j} (-1)^{\abs{J}}V_J,\quad j\geq 0. 
\endeq
It was shown in~\cite{Melin} that the expansion
\begeq
\label{Melin}
\varphi(H_0+V)\sim \sum_{j=0}^{\infty} \frac{ \varphi^{(j)}(H_0)V_j}{j!},\quad \varphi\in {\cal S}(\real^n), 
\endeq
holds in an appropriate space of pseudo-differential operators --- 
see~\cite{Melin} for the precise statement. 
\vskip 2mm
\noindent
The following observation is now crucial. 

\begin{prop}
\label{mc1}
The operators $V_j$ satisfy the following recurrent relation, 
\begin{equation}
\label{rec}
V_0=I; \quad V_j=V_{j-1}V+[V_{j-1},H_0],
\end{equation}
and are given explicitly by 
\begin{equation}
\label{vj1}
V_j=\sum_{k=0}^j (-1)^k {j\choose k} H_0^k H^{j-k}.
\end{equation}
We have that
\begeq
\label{decayvj}
V_j=\sum_{\abs{\nu}\leq j-1} p_{\nu}(x)\partial_x^{\nu},\quad p_{\nu}\in S^{-\varepsilon j}(\real^n). 
\endeq
\end{prop}
\begin{proof}
This result is very close to Theorem 2.3 in~\cite{Pol1}, and for the sake of completeness we present the proof, which proceeds by induction. 
Introduce the set
$$
W_j=\{J=(j_1,..,j_k):\norm{J}=j\}.
$$
Then $V_0=I$ by definition, and $V_1=V$ since $W_1$ consists of a single vector $(0)$. In order to pass from $V_{j-1}$ to $V_j$, we argue as 
follows. Let $W_j^0$ denote the set of all vectors in $W_j$ whose last entry is $0$. Consider a mapping $p:W_j \to W_{j-1}$ 
defined as follows: for $J=(j_1,..,j_{k-1},0) \in W_j^0$ we have $p(J)= (j_1,..,j_{k-1})$ and for $J=(j_1,..,j_k) \in W_j\setminus W_j^0$ 
we have $p(J)=(j_1,..,j_k-1)$. It is clear that the mapping $p$  is well--defined. Moreover, it is a surjection and every element of 
$W_{j-1}$ has exactly two preimages --- one in $W_j^0$ and one in its complement. Applying the definition of the operator $V_j$ we prove 
the inductive step --- $V_{j-1}V$ is the contribution of vectors from $W_j^0$ and $[V_{j-1},H_0]$ is the contribution of vectors belonging 
to $ W_j\setminus W_j^0$. This completes the proof of the first part of the theorem. 
For obtaining the closed formula (\ref{vj1}), let us rewrite the recurrent relations (\ref{rec}) in the following way:
$V_1=H_0-(H_0-V)=V$; $V_j=V_{j-1}H-H_0V_{j-1}$. Proceeding by induction and recalling that ${j \choose k}+{j \choose k-1}={j+1 \choose k}$ 
we get (\ref{vj1}). Finally, a simple proof by induction using (\ref{rec}) gives (\ref{decayvj}). 
\end{proof}

\Remark. In view of (\ref{Melin}), the expression (\ref{vj1}) should be compared with $V^j=\left((H_0+V)-H_0\right)^j$. 

\vskip 2mm
As in Section 1, set  
$N:=[\frac{n}{\varepsilon}]$, the integer part of $n/{\varepsilon}$. This value of $N$ will be kept fixed throughout this section. 
It follows from Theorem 2.4 in~\cite{Melin} that if 
$\varphi \in S^{0}(\real)$,i.e. if 
$(1+\abs{\lambda})^{k} \varphi^{(k)}(\lambda)$ is bounded on $\real$ for every $k$, then the operator 
\begeq
\label{mel1}
\varphi(H)-\sum_{m=0}^{N}\varphi^{(m)}(H_0)V_m/m!
\endeq
is of trace class. Its trace is equal to 
\begeq
\label{trace}
T_V(\varphi)+\int_0^{\infty} \varphi(\lambda) f(\lambda)\,d\lambda,
\endeq
where $T_V\in {\cal E}'(\real)$, and the singular support of $T_V$ is contained in the compact set $\Spec_{\mathrm{p}}(H)\cup\{0\}\subset 
\overline{\real_-}$. Here $\Spec_{\mathrm{p}}(H)$ is the point spectrum of $H$, and we also use the fact that 
$H$ has no positive eigenvalues. The trace distribution $f(\lambda)$ could be viewed as a generalization of the scattering phase, 
see also Proposition \ref{pr}. We have that $f(\lambda)\in S^{-1-\delta}$ for some $\delta>0$, and there is an asymptotic expansion of 
$f(\lambda)$ in a neighbourhood of infinity. To describe this expansion, we recall that when $V_m(x,\xi)$ is the full symbol of $V_m$, the 
following function was introduced in~\cite{Melin}, 
\begeq
\label{am}
{\rm a}_m(t)=2^{-1} (2\pi)^{-n} t^{(n-2)/2}\int\!\!\!\int V_m(x, t^{1/2}\omega)\,dx\,d\omega,\quad m>N.
\endeq
Here the integration is performed over $\real^n\times S^{n-1}$. We notice that the integral converges in view of (\ref{decayvj}). Set now
$$
b_m(t)=(-1)^m {\rm a}_m^{(m)}(t).
$$
It follows then from (\ref{am}) that, for each $m>N$, the following powers of $t$ occur in $b_m(t)$: $t^{n/2-1+\abs{\nu}/2-m}$, where 
$\abs{\nu}\leq m-1$ and all entries in $\nu$ are even. A combination of Theorem 2.4 in~\cite{Melin} with this observation shows that the 
following expansion is true, as $\lambda \rightarrow \infty$, 
\begeq
\label{exp}
f(\lambda)\sim \sum_{j=1}^{\infty} \beta_j \lambda^{n/2-1-j},\quad n\quad\wrtext{is odd}. 
\endeq
where 
\begeq
\label{killfirst}
\beta_j=0, \quad\wrtext{for}\quad j<(N+2)/2.
\endeq
In the case when $n$ is even, it follows that the coefficients $\beta_j$ are different from $0$ only for $j\leq n/2-1$. Taking 
(\ref{killfirst}) into account, and using that $n<N+4$, we get that 
\begeq
\label{expeven}
f(\lambda)={\cal O}(\lambda^{-\infty}),\quad n\:\:\wrtext{is even}.  
\endeq

\Remark. The estimate (\ref{expeven}) should be compared with the fact that, in the short-range case, it is true that when $n=2$,
$$
s'(\lambda)={\cal O}(\lambda^{-\infty}),\quad \lambda\rightarrow \infty. 
$$

\noindent
\Remark. We remark that trace distributions associated to Schr\"odinger operators with long range potentials were also recently considered 
in~\cite{Bouclet1},~\cite{Bouclet2}. The approach used there is different, and is based on~\cite{Koplienko}. In particular, the operators 
$V_j$ playing a major r\^ole 
in~\cite{Melin} and in the present paper, do not occur in~\cite{Bouclet1},
\cite{Bouclet2}. 

\vskip 2mm
We are going to derive concise formulas for the coefficients $\beta_j$ in (\ref{exp}), and, similarly to Theorem 3.1, this will be done by 
considering the trace class operator 
\begeq
\label{trreg}
e^{-tH}-\sum_{m=0}^{N} (-1)^m \frac{t^m}{m!}e^{-tH_0} V_m, \quad t>0.
\endeq
>From the results of~\cite{Melin} we know that this is a pseudo-differential operator with integrable symbol, so that its trace is given by  
\begeq
\int \left(e^{-tH}(x,x)-\sum_{m=0}^{N} (-1)^m \frac{t^m}{m!}\left(e^{-tH_0} V_m\right)(x,x)\right)\,dx. 
\endeq
We would like to rewrite the integrand so that it would involve the operators $X_m$ introduced in (\ref{xj}). 

\begin{lemma}
We have 
$$
\left(e^{-tH_0}V_m\right)(x,x)=(-1)^m \left(X_m e^{-tH_0}\right)(x,x). 
$$
\end{lemma}
\begin{proof}
Using (\ref{vj1}), by repeated partial integrations, we get that the kernel of $e^{-tH_0}V_m$ equals 
\begin{eqnarray*}
& & \sum_{k=0}^m (-1)^k {m\choose k} H^{m-k}_y H_{0,y}^k e^{-tH_0}(x,y)=(-1)^m  
\sum_{k=0}^m (-1)^k {m\choose k} H^{k}_y H_{0,y}^{m-k} e^{-tH_0}(x,y) \\
& & =(-1)^m X_{m,y} e^{-tH_0}(x,y).
\end{eqnarray*}
Specializing to the diagonal, we get the statement of the lemma. 
\end{proof}

In view of Lemma 4.2, the trace of (\ref{trreg}) is given by 
\begeq
\label{tracexj}
\int \left(e^{-tH}(x,x)-\sum_{m=0}^{N} \frac{t^m}{m!}\left(X_m e^{-tH_0}\right)(x,x)\right)\,dx. 
\endeq
When computing the coefficients in the expansion of (\ref{tracexj}), we may use Theorem \ref{hh}. 

We are now ready to prove Theorem \ref{regt}.

\vskip 2mm
\noindent 
{\em Proof of Theorem 1.2:}
We shall first show that the expansion
\begeq
\label{local1}
e^{-tH}(x,x)-\sum_{m=0}^N \frac{t^m}{m!} \left(X_m e^{-tH_0}\right)(x,x)\sim (4\pi t)^{-n/2} \sum_{j=1}^{\infty} \alpha_j(x) t^j,
\endeq
holds, as $t \rightarrow 0^+$, and that $\alpha_j(x)\in L^1$. 

In proving (\ref{local1}) we notice that (\ref{alpha1}) follows from 
Theorem \ref{hh}. Furthermore, from the proof of Theorem \ref{hh} 
we recall that 
the powers of $t$ which occur in the sum in the left-hand side of 
(\ref{local1}) are $m-\abs{\nu}/2-n/2\leq N-n/2$, since $m\leq N$ 
and $0\leq \abs{\nu}\leq m-1$. Therefore, for $j\geq N+1$, we have that 
$\alpha_j(x)=a_j(x)$, which is (\ref{alpha2}). It remains to 
consider the case when $(N+2)/2\le j\leq N$. Arguing as in the proof of 
Theorem \ref{hh} we see that for these values of $j$, the coefficient in 
front of $t^j$ in the sum in the left-hand side of (\ref{local1}) is equal to $(4 \pi t)^{-n/2}$ times 
\begeq
\label{term}
\sum_{m=j}^{N} \frac{1}{m!} \sum_{\abs{\mu}=m-j} X_{m,y}\left(\frac{(y-x)^{2\mu}}{(2\mu)!}\right)\bigg |_{y=x} 
\frac{(-1)^{\abs{\mu}}(2\mu)!}{4^{\abs{\mu}} \mu!}.
\endeq
Expressions of this kind were studied in Sections 3 and 4 of~\cite{Pol3}, and it was shown that the combinatorial coefficients there can be 
simplified. Applying the arguments of~\cite{Pol3} to (\ref{term}) we get that (\ref{term}) is equal to 
$$
(-1)^j\sum_{k=0}^{N-j}
{N-j+\frac{n}{2}\choose k+\frac{n}{2}} 
\frac{H_y^{k+j}\left(d(x,y)^{2k}\right)|_{y=x}}{4^k k! (k+j)!}, 
$$
and then (\ref{alpha3}) follows. It remains to show that $\alpha_j(x)\in L^1$. In doing so we recall from (\ref{eq3}) that $a_j$ is obtained 
through contributions of terms involving $X_m$, for $j\leq m\leq 2j-1$. It follows that $a_j(x)\in L^1$, for $j\geq N+1$, 
since arguing as in Proposition 4.1 we see that the coefficients in $X_m$ are in $S^{-\varepsilon m}(\real^n)$. 
When $(N+2)/2\leq j\leq N$, it follows from (\ref{eq3}) and (\ref{term}) that $\alpha_j$ is given by 
$$
\alpha_j(x)=\sum_{m=N+1}^{2j-1} \frac{1}{m!} \sum_{\abs{\mu}=m-j} X_{m,y}\left(\frac{(y-x)^{2\mu}}{(2\mu)!}\right)\bigg |_{y=x}
\frac{(-1)^{\abs{\mu}}(2\mu)!}{4^{\abs{\mu}} \mu!} \in L^1. 
$$
It follows that $\alpha_j(x)\in L^1$ for every $j$. 

An application of Lemma \ref{aux} below then shows that the asymptotic expansion (\ref{local1}) can be integrated over $\real^n$, and this completes the
proof. 
{\hfill$\Box$}
\begin{lemma}
\label{aux}
When $V\in S^{-\varepsilon}(\real^n)$, we have that 
\begeq
\label{integrable}
(4 \pi t)^{n/2} e^{-tH}(x,x)-\sum_{j=0}^M a_j(x)t^j=R_M(x,t),\quad M\geq N, 
\endeq
where 
\begeq
\label{remainder}
\abs{R_M(x,t)}\leq h_M(x)t^{M+1},\quad h_M(x)\in L^1(\real^n). \quad 0\leq t \leq 1. 
\endeq
\end{lemma}
The proof of Lemma \ref{aux} is given in Section 6. 

\vskip 2mm
\noindent
We illustrate Theorem \ref{regt} by the following example. 

\smallskip

\noindent{\it Example.}
Consider the Schr\"odinger operator $H=-D_x^2+V(x)$ on the real line, where 
\begeq
\label{pot}
V(x)=\frac{1}{(1+x^2)^{1/6}}.
\endeq
We have that $V\in S^{-1/3}(\real)$, and therefore we take $N=3$. Explicit computations using Theorem \ref{regt} show that the coefficients 
$\alpha_j(x)$ are given by
$$\alpha_1(x)=0, \quad \alpha_2(x)=0,\quad
\alpha_3(x)=-\frac{1}{4}V'\,^2 - \frac{1}{3}V V'' + \frac{3}{20} V^{(4)},$$
and $\alpha_j(x)=a_j(x)\in L^1(\real)$ for $j\ge 4$. Now a computation shows that the first three local heat invariants are given by
$$a_1(x)=-V, \quad a_2(x)=\frac{1}{2} V^2-\frac{1}{6}V'', \quad
a_3(x)=-\frac{1}{6}(V^3 - \frac{1}{2}V'\,^2 - V V'' +\frac{1}{10} V^{(4)}),$$
and in view of (\ref{pot}), these functions are not integrable. In $a_3(x)$, the only non-integrable part is the first term involving $V^3$, 
which does not appear in $\alpha_3(x)$ while the other terms are present up to numerical factors. 

\vskip 2mm
We are now ready to compute the coefficients $\beta_j$ in the asymptotic expansion of 
the trace distribution $f(\lambda)$ defined in (\ref{trace}). 
We recall here that 
\begeq
\label{BKheat1}
\Tr\left(e^{-tH}-\sum_{j=0}^{N} (-1)^m \frac{t^m}{m!}e^{-tH_0}V_m\right)=\int_0^{\infty}e^{-t\lambda}f(\lambda)\,d\lambda+
T_V(\varphi_t),\quad t>0,
\endeq
where $T_V$ is a compactly supported distribution and $\varphi_t(s)=e^{-ts}$. 
\begin{prop}
\label{pr}
The following asymptotic expansion holds for 
$f(\lambda)$ as $\lambda \to \infty$,
\begeq
\label{odd1}
f(\lambda )\sim \sum_{j=1}^{\infty} \beta_j \lambda ^{n/2-j-1},\quad \lambda\rightarrow \infty,\;\;n\;\;\wrtext{odd}, 
\endeq
where
$$
\beta_j=(4 \pi)^{-n/2}\frac{\alpha_j}{\Gamma(n/2-j)}, \quad j=1,2,\ldots,\;\;
$$
and $\alpha_j$ are defined in  Theorem \ref{regt}.
We also have that 
\begeq
\label{even1}
f(\lambda)={\cal O}(\lambda^{-\infty}), \quad \lambda\rightarrow \infty,\;\;n\;\;\wrtext{even}. 
\endeq
\end{prop}
\begin{proof}
The result follows by an application of Theorem \ref{regt} and Lemma 5.2 in~\cite{CdV} to (\ref{BKheat1}), since, as 
$t\rightarrow 0$, $T_V(\varphi_t)$ has an asymptotic expansion in the 
integer powers of $t$.  
\end{proof}

We notice that Proposition \ref{pr} provides a long range version of the formula (3.7), which emphasizes the analogy between the trace 
distribution $f(\lambda)$ and the derivative of the scattering phase in the short range case.

\section{Taylor expansions for the heat semigroup}
\setcounter{equation}{0}
\setcounter{dref}{0}
Computations of heat invariants in~[15-17] and in 
the present paper were based on the commutator expansions of~\cite{AK} and 
~\cite{Melin}. A crucial observation was that the multiple commutators, present in these expansions, could be 
eliminated, and, in particular, this led to Theorem \ref{hh}. In this section we wish to point out that an expansion closely related to that of Theorem \ref{hh} 
can be derived from the usual Taylor formula for the semigroups $e^{-tH}$ and $e^{-tH_0}$. When stating the result, it seems instructive to proceed in 
the general context of semigroup theory. 

Let therefore $A$ and $B$ be two unbounded self-adjoint operators acting in a Hilbert space ${\cal H}$, generating strongly continuous 
semigroups $e^{tA}$ and $e^{tB}$, $t\geq 0$. When $m\geq 0$, introduce the operator 
\begeq
\label{Cm}
C_m(A,B)=\sum_{k=0}^m {m \choose k} A^k (-B)^{m-k} 
\endeq
with the domain 
$$
{\cal D}(C_m(A,B))=\cap_{k=0}^m {\cal D}(A^k B^{m-k}).
$$
The space of $C^{\infty}$-vectors for $A$ is denoted ${\cal D}^{\infty}(A)$, and similarly for $B$,
$$
{\cal D}^{\infty}(A)=\cap_{n=1}^{\infty} {\cal D}(A^n).
$$
This is a dense linear subspace of ${\cal H}$. 

\begin{theo} 
Let $A$ and $B$ be generators of strongly continuous semigroups such that ${\cal D}^{\infty}(B)\subset {\cal D}^{\infty}(A)$.
Then for any positive integer $N$ we have 
\begeq
\label{kant}
e^{tB}\varphi=\sum_{m=0}^N (-1)^m \frac{t^m}{m!} e^{tA} C_m(A,B)\varphi+R_N(t;A,B)\varphi, 
\endeq
for all $\varphi\in {\cal D}^{\infty}(B)$. Here the remainder $R_N(t;A,B)\varphi={\cal O}(t^{N+1})$ in ${\cal H}$.  
\end{theo}

\Remark. Specializing to the case when $B=-H$ and $A=-H_0$ in $L^2(\real^n)$, we see that $C_m(A,B)$ are precisely the operators $V_m$ introduced in 
Section 4. Theorem 5.1 should then be compared with the following expansion which follows from Theorem \ref{hh} and Lemma 4.2, 
$$
e^{-tH}(x,x)\sim \sum_{m=0}^{\infty} (-1)^m \frac{t^m}{m!} \left(e^{-tH_0}V_m\right)(x,x). 
$$

\begin{proof}
When $\varphi\in {\cal D}^{\infty}(B)\subset {\cal D}^{\infty}(A)$, it is true that $e^{tB}\varphi$ and $e^{tA}\varphi$ are 
$C^{\infty}$-functions of $t\in [0,\infty)$ with values in ${\cal H}$, with Taylor expansions at $t=0$, 
\begeq
\label{Taylor}
e^{tB} \varphi \sim \sum_{l=0}^{\infty} \frac{(tB)^l}{l!} \varphi,\quad e^{tA}\varphi \sim \sum_{l=0}^{\infty} \frac{(tA)^l}{l!}\varphi.
\endeq

If we substitute (\ref{Cm}) in the sum in the right-hand side of (\ref{kant}) we get that it is equal to 
\begeq
\label{heat3}
\sum_{m=0}^{N} \frac{t^m}{m!}e^{tA} \sum_{k=0}^m {m\choose k} (-A)^k (B)^{m-k} \varphi=
e^{tA} \sum_{m=0}^N \frac{(-tA)^m}{m!} \sum_{m=0}^N \frac{(tB)^m}{m!}\varphi+ R_1(t;A,B)\varphi,
\endeq
where
$$
R_1(t;A,B)\varphi=-\sum_{m=N+1}^{2N} t^m e^{tA} \sum_{k=m-N}^N \frac{(-A)^k}{k!} \frac{B^{m-k}}{(m-k)!}\varphi={\cal O}(t^{N+1})
\;\; \wrtext{in}\;\; {\cal H}. 
$$
Now, using (\ref{Taylor}) we see that 
$$
e^{tA} \sum_{m=0}^N \frac{(-tA)^m}{m!} \sum_{m=0}^N \frac{(tB)^m}{m!}\varphi=e^{tB}\varphi+R_2 (t;A,B)\varphi,
$$
where $R_2(t;A,B)={\cal O}(t^{N+1})$ in ${\cal H}$. We have therefore proved that 
$$
e^{tB}\varphi-\sum_{m=0}^N (-1)^m \frac{t^m}{m!} e^{tA} C_m(A,B)\varphi ={\cal O}(t^{N+1}) \;\;\wrtext{in}\;\;{\cal H}.
$$
\end{proof}

\Remark. A somewhat different proof of Theorem 5.1 leading to a more explicit expression for the remainder in (\ref{kant}), was given 
in~\cite{Kantor1}, see also~\cite{Kantor2}. The expansion (\ref{kant}) was referred to in~\cite{Kantor1} as a non-commutative Taylor formula.
  
\section{Proof of Lemma \ref{aux}}  
\setcounter{equation}{0}
\setcounter{dref}{0}
In this section we shall prove Lemma 4.3. 
We may assume that $M=N$. When proving (\ref{remainder}), we shall make use of the Feynman-Kac formula (\ref{eq2}), which we rewrite in the 
following form, see~\cite{Simon}, 
\begeq
\label{feynmankac}
e^{-tH}(x,x)=(4 \pi t)^{-n/2} E\left(\exp\left(-t\int_0^1 V(x+\sqrt{2t}b(s))\,ds\right)\right). 
\endeq
Here $b(s)$, $0\leq s\leq 1$, is the Brownian bridge which returns to $0$ at time $t=1$. Taking the Taylor 
expansion of the exponential in (\ref{feynmankac}) we get 
\begeq
\label{expand}
(4 \pi t)^{n/2} e^{-tH}(x,x)=\sum_{k=0}^N \frac{(-t)^k}{k!} E\left(\left(\int_0^1 V(x+\sqrt{2t}b(s))\,ds\right)^k\right)+
R_N(t,x),
\endeq
where 
\begeq
\label{rest}
\abs{R_N(t,x)}\leq \frac{e^{\norm{V}_{L^{\infty}}}}{(N+1)!} t^{N+1} E\left(\left(\int_0^1 \abs{V(x+\sqrt{2t}b(s))}\,ds\right)^{N+1}\right), 
\quad 0\leq t\leq 1.
\endeq
An elementary estimate of the expectation in (\ref{rest}) shows that it is bounded by a constant times $(1+\abs{x})^{-\varepsilon(N+1)}$, which is 
integrable, since $N+1>n/{\varepsilon}$. Consider now the expression 
\begeq
\label{expect}
E\left(\left(\int_0^1 V(x+\sqrt{2t}b(s))\,ds\right)^k\right),\quad 1\leq k\leq N.
\endeq
Since $E(b(s)^{\nu})=0$ unless all indices in $\nu=(\nu_1,\ldots \nu_n)$ are even, this has a complete asymptotic expansion in 
integer powers of $t$, and we are interested in showing that (\ref{expect}) has an expansion of the form 
\begeq
\label{want}
\sum_{j=0}^{N-k} a_{k,j}(x)t^j +R_{N,k}(x,t),
\endeq
for some $a_{k,j}(x)$, where
$$
\abs{R_{N,k}(x,t)}\leq h_{N,k}(x)t^{N+1-k},\quad h_{N,k}\in L^1.
$$
We shall apply Taylor's formula 
\begin{eqnarray}
& & V(x+\sqrt{2t}b(s))=\sum_{\abs{\nu}\leq 2N} \frac{\partial^{\nu} V(x)}
{\nu!}
(2t)^{\abs{\nu}/2} b(s)^{\nu} \\ \nonumber 
& + & (2N+1)\sum_{\abs{\nu}=2N+1} 
\frac{(2t)^{\abs{\nu}/2} b(s)^{\nu}}{\nu!} \int_0^1 \partial^{\nu}V(x+y\sqrt{2t}b(s))(1-y)^{2N}\,dy, 
\end{eqnarray}
which we rewrite in the following way, 
\begeq
\label{Taylor1}
V(x+\sqrt{2t}b(s))=\sum_{j=0}^{2N} f_j(x,s)t^{j/2}+t^{N+1/2}R_N(x,t,s).
\endeq
We have that $f_j(\cdot,s)$ and $R_N(\cdot,t,s)$ are ${\cal O}(1)(1+\abs{x})^{-\varepsilon-j}$ and ${\cal O}(1)(1+\abs{x})^{-\varepsilon-2N-1}
\in L^1$, respectively. We substitute (\ref{Taylor1}) in (\ref{expect}) and expand the $k$-th power there by means of the multinomial theorem. It
 follows that 
(\ref{expect}) is equal to a linear combination of terms of the form 
\begeq
\label{multi}
E \biggl((g_0(x))^{k_0}(g_1(x))^{k_1} t^{k_1/2}\ldots t^{N k_{2N}}(g_{2N}(x))^{k_{2N}}
t^{k_{2N+1}(N+1/2)}\left(\int_0^1 R_N\,ds\right)^{k_{2N+1}}\biggr).
\endeq
Here $k_0+k_1+\ldots k_{2N+1}=k$, and we have put 
$$
g_j(x)=\int_0^1 f_j(x,s)\,ds, \quad j=0,\ldots 2N.
$$
If $k_{2N+1}\geq 1$, then the term (\ref{multi}) has the desired bound of the form $t^{N+1-k}$ times an integrable functi
on of $x$. 
When analyzing the expression (\ref{multi}) we may therefore assume that $k_{2N+1}=0$. We then get terms of the form 
\begeq
\label{multi2}
t^{(1/2)\sum_{j=0}^{2N}j k_j}E\biggl(g_0(x)^{k_0}\ldots g_{2N}(x)^{k_{2N}}\biggr),
\endeq
with $\sum_{j=0}^{2N} k_j=k$. The terms which contribute to the remainder in (\ref{want}) are those for which 
$$
\half \sum_{j=0}^{2N}j k_j \geq N+1-k,
$$
and the coefficient in (\ref{multi2}) in front of such a power of $t$ can be estimated by a constant times $(1+\abs{x})^{-1}$ raised to the power
\begeq
\label{decay}
\varepsilon k+\sum_{j=0}^{2N} j k_j \geq \varepsilon k+2(N+1)-2k \geq 2+\varepsilon N >n+2-\varepsilon>n,
\endeq
which therefore gives an integrable bound. In (\ref{decay}) we have used that $k\leq N$. This completes the proof of Lemma \ref{aux}. 
{\hfill$\Box$}

\end{document}